\documentclass[a4paper,12pt]{article}
\usepackage{amsmath,amssymb,amsfonts,amsthm}
\usepackage[english]{babel}
\usepackage{color}
\usepackage{graphicx}
\usepackage{wrapfig}
\usepackage{rcdthm}
\usepackage{rcdcom}
\usepackage[top=2.5cm,left=3cm,right=3cm,bottom=2.5cm]{geometry}

\date{}

\begin{document}
\newcommand{\beq}{\begin{equation}}
\newcommand{\eeq}{\end{equation}}
\newcommand{\bOm}{\bs\Omega}
\newcommand{\bom}{\bs\omega}
\newrcdthm[example]{exl}{Example}
\newrcdthm{pro}{Proposition}

\title{\bf Dynamics and Control of a~Spherical Robot with an Axisymmetric Pendulum Actuator\vspace{1cm}}
\author{Tatyana\,B.\,Ivanova\footnote{tbesp@rcd.ru},\, Elena\,N.\,Pivovarova\footnote{archive@rcd.ru}\vspace{5mm}\\
\small Udmurt State University\\
\small ul. Universitetskaya 1, Izhevsk, 426034 Russia}

\maketitle

\begin{abstract}
This paper investigates the possibility of the motion control of a~ball
with a~pendulum mechanism with non-holonomic constraints using gaits~---
the simplest motions such as  acceleration and deceleration during the
motion in a~straight line, rotation through a~given angle and their
combination. Also, the controlled motion of the system along a~straight
line with a~constant acceleration is considered. For this problem the
algorithm for calculating the control torques is given and it is shown
that the resulting reduced system has the first integral of motion.
\end{abstract}

\section*{Introduction\label{vved}}
\addcontentsline{toc}{section}{Introduction}

For the last ten years an enormous amount of research has been devoted to
the dynamics and control of such vehicles as single wheel robot and
spherical robot moving due to  changes in the position of the center of
mass (see, for example,
\cite{komarov,mart,mart2,das,Kayacan,Michaud,nagai,Schroll}). Also, the
possibility of controlling such systems by using other internal
mechanisms, for example, rotors~\cite{bkm1,bkm}, is actively studied. The
interest in such systems is determined by the presence of some advantages
in maneuverability over wheeled vehicles. For a~detailed literature review
on spherical robots with various moving mechanisms, their description and
application areas see~\cite{chase,nagai,Schroll}.

The motion of spherical robots moving due to pendulum oscillations is
studied in~\cite{bm,komarov,mart2,pi,Kayacan,Michaud,nagai,Schroll}. In
particular, Nagai~\cite{nagai} considers the control of the motion of a~pendulum spherical robot on an inclined plane and finds the maximum
inclination angle of the plane at which the vehicle can move up the plane
(a~similar problem is considered in~\cite{mart2} for a~single wheel
robot). Schroll~\cite{Schroll} addresses the problem of obstacle
negotiation and finds the maximum height of the obstacle  which such a~robot  can overcome. In~\cite{Kayacan} Kayacan et al. consider the control
of the motion of a~ball in a~straight line and a~circle for various types
of controllers.

In~\cite{komarov} the authors consider the problems of controlling a
spherical robot with a~pendulum drive in the case of slipping at the
contact point and in the non-holonomic case. The focus is on determining
an algorithm for an additional control to approach the given trajectory
from an arbitrary point in the case of a singular matrix defining the
control. However, this investigation~\cite{komarov} leads to a physically
strange conclusion that introduction of a small parameter can lead to the
controllability of this system, whereas vanishing of this parameter
corresponds to the absence of controllability, which casts some doubts on
the accuracy of the results obtained.


This work is a~continuation of  analysis of the dynamics of a~spherical
shell rolling without slipping on a~horizontal plane with Lagrange's top
fixed at the center of the shell~\cite{bm,pi}. Previously in~\cite{bm} the
equations of motion for a~free system were obtained (the equations of
motion for the system situated inside the rolling ball were also obtained
in~\cite{bolotin}), all necessary first integrals and an invariant measure
were found, the reduction to quadratures is given. In~\cite{pi} the
\textit{free} motion of the ball with Lagrange's top was considered, the
stability of  periodic solutions was analyzed and the trajectories of the
contact point were constructed.

This paper investigates the \textit{controlled} motion of a~ball with a~pendulum. Special attention is paid to  the control of the ball using
gaits --- the simplest motions (such as  acceleration to a~certain
velocity and deceleration in a~straight line and  rotation through a~given
angle) and their combination for one oscillation of the pendulum and
specific examples of such motions are shown. In the last part of the paper
we consider the controlled motion of the ball along a~straight line with
fixed parameters according to a~predetermined law of motion, indicate
shortcomings of this approach and consider an algorithm for computing the
control torques by example.

\section{Equations of motion}

\begin{wrapfigure}[13]{o}{65mm}
\vspace{-4mm}
\begin{center}
\includegraphics[scale=0.9]{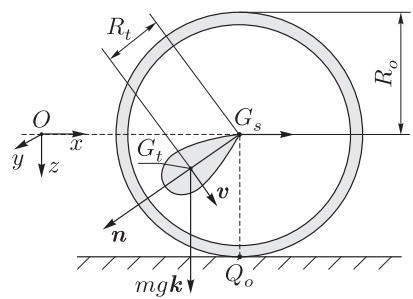}
\vspace{-0.5cm}
\parbox{55mm}{\caption{Spherical shell with an axisymmetric pendulum fixed at its center.%
}}\end{center}
\label{ivanova-fig1}
\end{wrapfigure}We consider a~spherical shell (Fig.~\ref{ivanova-fig1}) relative to a~fixed reference frame (the axis $Oz$ is directed vertically downwards).
Let $G_s$ be the center of mass of the shell, $G_t$ the center of mass of the top and let $R_t=|G_sG_t|$ denote the distance between them. The vector $\bs n$ is directed along the axis of symmetry of the pendulum (here and in the sequel, the vectors are denoted
by  bold italic letters).

Assuming that in the system of the top’s
principal axes its tensor of inertia is ${\bf \hat{i}}=\diag(i,i,i+j)$, we can represent the kinetic energy
of the system as~\cite{bm,pi}:
\begin{equation*}
T=\dfrac12\left(M\bV^2+I\bs\Omega^2\right)+\dfrac12\left(m\bs v^2+i\bs\omega^2+j(\bs\omega,\bn)^2\right)\!,
\end{equation*}
where $\bV$ and $\bs\Omega$ are the linear velocity of the center of the shell and the shell’s angular velocity, $M$
and $I$ stand for its mass and moment of inertia, $\bs v$ and $\bs\omega$ are the linear velocity of the center of
mass of the top and its angular velocity, and  $m$ is the mass of the top.

In this paper we will investigate the controlled motion of the system. The control torque  can be generated by the engine which is installed at the point of attachment of the pendulum to the ball and sets the pendulum and, accordingly, the ball in motion. Let $\bs Q$ be the torque generated by the engine.

Then the change of the angular momentum relative to $G_s$ and the momentum of the shell can be written
as\nopagebreak\vspace{1mm}
\beq
\label{seq1}
\begin{aligned}
\dfrac{d}{dt}\dfrac{\partial T}{\partial\bs\Omega}&=I\dot{\bs\Omega}=R_o\bk\times\bs{N}_o-\bs
Q,\\[2mm]
\dfrac{d}{dt}\dfrac{\partial T}{\partial\bV}&=M\dot{\bV}=\bs{N}_o+\bs{N}_t+Mg\bk,
\end{aligned}
\eeq
where $\bs{N}_o$ and $\bs{N}_t$ are the reaction forces acting on the shell at the
contact point $Q_o$ and the point  of attachment of the top $G_s$. For the top relative
to its center of mass we have\nopagebreak\vspace{1mm}
\beq
\label{seq2}
\begin{aligned}
\dfrac{d}{dt}\dfrac{\partial
T}{\partial\bs\omega}&=i\dot{\bs\omega}+j(\bs\omega,\bn)\dot\bn=R_t\bn\times\bs{N}_t+\bs
Q,\\[2mm]
\dfrac{d}{dt}\dfrac{\partial T}{\partial\bs v}&=m\dot{\bs v}=mg\bk-\bs{N}_t,
\end{aligned}\vspace{2mm}
\eeq
The linear velocity of the center of mass of the top is determined by the relation
\[
\bs v=\bV+R_t\bs\omega\times\bn=R_o\bk\times\bs\Omega+R_t\bs\omega\times\bn.
\label{v}
\]
The linear velocity of the center of mass of the shell $\bs V$ is related with the angular velocity of the shell $\bOm$ by the no-slip condition at the contact point $Q_o$:
\beq
\bV=\dot{\br}_s=R_o\bk\times\bs\Omega,
\label{rad-vector}
\eeq
where $\bs r_s$ is the radius vector of the contact point.

The evolution of the vector $\bn$ can be found from the equation
\begin{equation*}
\dot{\bn}=\bs\omega\times\bn.
\end{equation*}

Eliminating the reaction forces ${\bN}_o$ and ${\bN}_t$ from Eqs.~\eqref{seq1} and \eqref{seq2}, we obtain the equations of controlled motion for the spherical shell with the axisymmetric top fixed at its geometrical center:
\beq
\label{eqmQ}
\begin{aligned}
&{\bf J}\dot{\bs\Omega}+mR_o^2\bk\times(\dot{\bs\Omega}\times\bk)-mR_oR_t\bk\times(\dot{\bs\omega}\times\bn)=mR_oR_t\bk\times(\bs\omega\times\dot{\bn})-\bs
Q,\\[2mm]
&i\dot{\bs\omega}+mR_t^2\bn\times(\dot{\bs\omega}\times\bn)-mR_oR_t\bn\times(\dot{\bs\Omega}\times\bk)={}\\[1mm]
&\hspace*{3.4cm}=-j(\bs\omega,\bn)\dot{\bn}-mR_t^2\bn\times(\bs\omega\times{\dot \bn})+mgR_t\bn\times\bk+\bs Q,\\
&\dot{\bn}=\bs\omega\times\bn,
\end{aligned}
\eeq
where ${\bf J}=\diag(I+MR_o^2,I+MR_o^2,I)$, $\bk=(0,0,1)^{\rm T}$.\goodbreak

Differentiating Eq.~\eqref{rad-vector} with
respect to time, we obtain the expression for $\dot{\bs\Omega}$:
\beq
\dot{\bOm}=\dfrac{1}{R_o}\,\dot{\bs V}\times\bk=\dfrac{1}{R_o}\,\bs a\times\bk,
\label{dOmega}
\eeq
where $\bs a$ is the acceleration of the contact point.

Our goal is to determine a~control torque $\bs Q$ such that the contact
point (and thus the center of the ball) moves according to the specified
law of motion $\bs r_s(t)=(x(t),y(t),0)^{\rm T}$. The velocity of the
center of mass $\bV(t)=(\dot{x}(t),\dot{y}(t),0)^{\rm T}$ and its
acceleration $\bs a(t)=(\ddot{x}(t),\ddot{y}(t),0)^{\rm T}$ are prescribed
functions as well.

By prescribing the law of motion of the contact point on a~plane $\bs
r_s(t)$ from~\eqref{eqmQ} and \eqref{dOmega}, we obtain a~system of nine
equations for the projections of the vectors~$\bn$,~$\bom$,

\noindent
and $\bs Q$
which can be represented as \beq {\bf F}\bs y={\bs A}, \label{matrF} \eeq
where $\bs y = (\dot{\bs\omega},\dot{\bn},\bs Q)^{\rm T}$, ${\bf F}$ is
the matrix whose elements depend on $\bs\omega$ and $\bn$; $\bs A$ is the
vector function depending on $\bs\omega,\bn$, and $\dot{\bOm}$.

For the system~\eqref{matrF} to have a~solution, i.e., for a~controlled
motion to be possible, it is necessary that there exist the inverse matrix
${\bf F}^{-1}$, i.e., that the condition ${\det{\bf F}\neq0}$ be satisfied,
which requires the fulfilment of the condition (see also~\cite{komarov})
\beq i+mR_t^2>mR_oR_t. \label{nerav} \eeq

The condition~\eqref{nerav} can be satisfied by choosing the corresponding geometrical characteristics of the system. For definiteness  we choose as the pendulum a~thin disk suspended on a~massless rod, whose parameters satisfy the condition~\eqref{nerav}: the radius of the disk is $R_d=0.92R_o$ and the length of the rod is $R_t=0.25R_o$. Note that this system is similar to a~vehicle rolling without friction on the inner surface of a~spherical shell.

Thus, to determine the control torque $\bs Q(t)$, it is necessary to solve
the system~\eqref{matrF} of differential equations with the corresponding
initial conditions. However, by\linebreak prescribing only the law of motion $\bs
r_s(t)$, the system cannot start a~new maneuver, for example, for changing
the motion direction when bypassing an obstacle. In addition, it is
necessary that the ball stop at the end point of a~trajectory, and this
requires that the pendulum be in the lower position and the velocity and
acceleration be equal to zero at the final instant of time. As a~rule, it
is very difficult to satisfy such conditions beforehand,  since this
substantially complicates the determination of the function $\bs r_s(t)$
for the entire trajectory of motion.

Another approach to the control in maneuvering along a~general trajectory is to use gaits. This approach implies that each motion must be performed \textit{for one oscillation of the pendulum}, which is the necessary condition for a~new motion to start. By combining such motions we can get any complex trajectory (which is useful, for example, for bypassing an obstacle).

In this paper we will consider both  approaches --- the control using gaits\linebreak (Section~\ref{geit}) and the motion control  with fixed parameters (such as acceleration, Section~\ref{fix}).

For ease of use, we will write all the equations of motion  in
dimensionless form. For this we take the mass of the pendulum $m$ as the
unit of mass, the quantity $t_0=\sqrt{\frac{i+j}{mgR_t}}$ as the unit of
time and the quantity $x_0=gt_0^2$ as the unit of length, where $g$ is the
free-fall acceleration, i.e., in the equations of motion we make the
changes
\[
\dfrac{t}{t_0}\rightarrow t,\quad\dfrac{x}{x_0}\rightarrow
x,\quad \dfrac{M}{m}\rightarrow M,\quad \bOm\,t_0\rightarrow \bOm,\quad
\bom t_0\rightarrow \bom, \quad \dfrac{\bs Q\,t_0^2}{mx_0^2}\rightarrow
\bs Q.
\]

 Also, to abbreviate some of our forthcoming formulae, we introduce the following notation:
\[
 i_0=i+R_t^2,\quad
 I_0=I+(1+M)R_o^2.
\]

\section{The control using gaits and their connection}\label{geit}

\begin{wrapfigure}[11]{o}{50mm}
\vspace{-4mm}
\centering\includegraphics[scale=0.9]{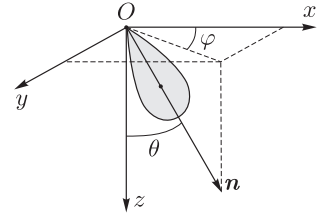}
\parbox{40mm}{\caption{\label{fig_theta}Determination of the angles $\theta$ and $\varphi$.}}
\end{wrapfigure}
We represent the vector $\bn$ directed along the axis of symmetry of the pendulum as
 \[
\bn=\left( \sin\theta\cos\varphi,\sin\theta\sin\varphi, \cos\theta\right)^{\rm T} ,
 \]
 where $\theta$ is the angle of deviation of the pendulum from the vertical and $\varphi$ is the angle between the axis~$Ox$ and the direction of oscillation of the pendulum (Fig.~\ref{fig_theta}).

 To determine $\bs Q(t)$, we will specify the angle of deviation of the pendulum from the vertical~$\theta$ such that at the start and the end of a~maneuver the pendulum is in the lower position:
 \beq\theta(t,\alpha,T)=\alpha\sin^2\left(\dfrac{\pi t}{T}\right)\!,\qquad
\left.\theta(\alpha,T)\right|_{t=0}=\left.\theta(\alpha,T)\right|_{t=T}=0.\label{theta}\eeq
 where $\alpha$ is an as yet unknown parameter defining the amplitude of oscillation,  $T$ is the specified time of one oscillation of the pendulum which is equal to the time of one maneuver.

To determine the parameter $\alpha$ corresponding to the specified change of velocity, it is necessary to express from the system~\eqref{eqmQ} using~\eqref{dOmega} and~\eqref{theta} the acceleration $a(t,\alpha,T)$ on which  the additional condition is imposed
\beq\Delta V(\alpha,T)=\int\limits_0^T a(t,\alpha,T)\,dt.\label{DeltaV}\eeq

Integrating~\eqref{DeltaV} for different values of the parameter $\alpha$ and the period of oscillations  $T$, we  obtain the dependence $\Delta V(\alpha,T)$ (quadric surface), from which we  find the parameter $\alpha$ by choosing the required velocity and the duration of the maneuver $T$.

Then, knowing the value $\alpha$, we can explicitly determine the acceleration $a(t,\alpha,T)$ and express the control torques from equations~\eqref{eqmQ}.

We will consider the described algorithm in the specific cases:  acceleration of the ball to a~given velocity in a~straight line and rotation for one oscillation of the pendulum.

\subsection{\label{razgon_po_pryamoi}Acceleration in a~straight line}

We direct the axis $Ox$ along the direction of motion. It is obvious that during the  motion in a~straight line the pendulum will oscillate in one plane, in this case in the plane $Oxz$, therefore, the vector $\bn$ directed along the axis of symmetry of the pendulum can be written as
\[
\bn=\left(
\sin\theta,
0,
\cos\theta
\right)^{\rm T},
\]
where $\theta$ is the angle of deviation of the pendulum from the vertical specified by the
expression~\eqref{theta}.

Only three of nine equations of motion~\eqref{eqmQ} are nontrivial:
\beq
\begin{aligned}
&I_0\dot{\Omega}_2-R_oR_t(\dot{\omega}_2\cos\theta-\omega_2^2\sin\theta)+Q_2=0,\\
&i_0\dot{\omega}_2-R_oR_t\dot{\Omega}_2\cos\theta+R_t\sin\theta-Q_2=0,\\
&\dot{\theta}=\omega_2.
\end{aligned}
\label{eq_acc}
\eeq

The value of $\dot{\Omega}_2$ is determined from~\eqref{dOmega}:
\beq
\dot{\Omega}_2=-\dfrac{a(t)}{R_o},\quad\dot{\Omega}_1=\dot{\Omega}_3=0
\label{dOm2}
\eeq

Substituting~\eqref{dOm2} into~\eqref{eq_acc}, we obtain the expression for acceleration of the ball in the form
\beq
a(t,\alpha,T)=R_o\dfrac{\ddot{\theta}(R_oR_t\cos\theta-i_0)+R_oR_t\dot{\theta}^2\sin\theta-R_t\sin\theta}{I_0-R_oR_t\cos\theta},\label{acc}\eeq
where $\theta$, $\dot{\theta}$ and $\ddot{\theta}$ are explicit functions of time and the parameters $\alpha$ and $T$ defined from~\eqref{theta}.

\begin{figure}[!b]
\begin{tabular}{cc}\small
\includegraphics{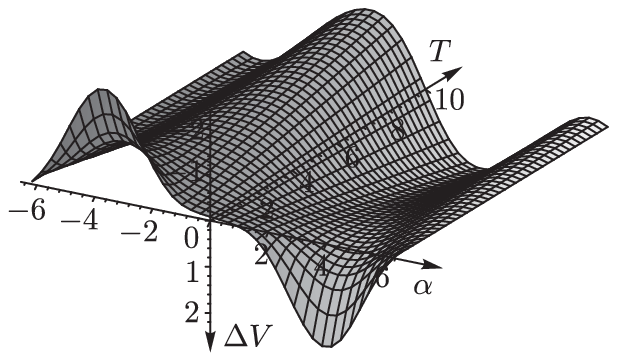}&
\includegraphics{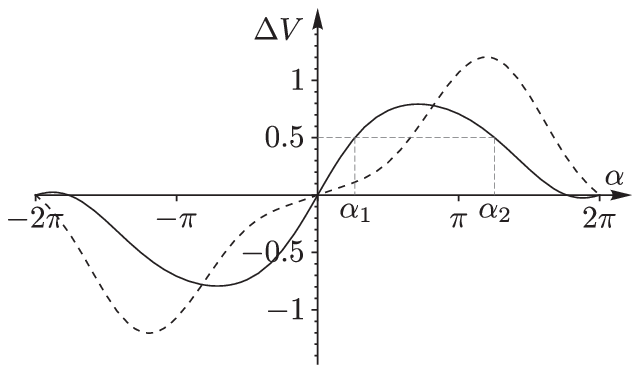}\\
({\it a})&({\it b})
\end{tabular}
\caption{\label{dv}({\it a}) The surface  $\Delta V(\alpha,T)$.  ({\it b}) Section
formed by the intersection of the surface $\Delta V(\alpha,T)$ with the
planes $T=5$ (solid line) and $T=1$ (dashed line).}
\end{figure}

Numerically integrating~\eqref{DeltaV} using~\eqref{acc}, we obtain the surface illustrated in Fig.~\ref{dv}({\it a}). This surface $\Delta V(\alpha,T)$ is antisymmetric relative to the change $\alpha\rightarrow-\alpha$ (this can be easily shown by substituting the function $\theta(t,\alpha_1,T)$  into~\eqref{acc} in the explicit form~\eqref{theta}).

In addition, $\Delta V_{\max}\rightarrow\infty$ as $T\rightarrow0$, i.e., the quicker the oscillation of the pendulum occurs, the more the velocity increases.

 Figure~\ref{dv}({\it b}) shows sections of the surface for different values of
 $T$. From these graphs we can see that at least two values of $\alpha$
 correspond to given $\Delta V$ and $T$: at~$\alpha=\alpha_1$ the pendulum
 performs an oscillation by a~smaller angle than at $\alpha=\alpha_2$
 (because of this the velocity varies non-monotonically at
 $\alpha=\alpha_2$).

Choosing the specific value of $\alpha$, for example $\alpha_1$,  substituting~\eqref{acc} into the equations of motion~\eqref{eq_acc} and numerically integrating them, we obtain the control torques in the form
\begin{gather*}
Q_1=Q_3=0,\\
Q_2=\dfrac{a(t,\alpha,T)}{R_o}\,I_0+R_oR_t(\cos{\theta}\ddot{\theta}-\sin{\theta}\dot{\theta}^2).
\end{gather*}

Thus, for acceleration in a~straight line it is necessary to generate a~control torque which is orthogonal to the direction of motion and to the plane of oscillations of the pendulum.

\begin{exl}
 We consider the acceleration of the ball from rest to the velocity $V_1=0.5$ for the time interval $T=5$. The value of the parameter $\alpha=0.83$ corresponds to such a~motion (see Fig.~\ref{dv}({\it b})). Since the ball moves from rest and the pendulum is in the lower position, the initial conditions have the form
 \[
\bn=\left(0,0,1\right)^{\rm T},\quad
\bs\omega=\left(0,0,0\right)^{\rm T},\quad
\bs r=\left(0,0,0\right)^{\rm T},\quad
\bs V=\left(0,0,0\right)^{\rm T}.
\]

Substituting the computed value $\alpha$ into the equations of
motion~\eqref{eq_acc} and\linebreak numerically integrating them (with the
specified initial conditions), we obtain the corresponding control $\bs
Q=(0,Q_2,0)^{\rm T}$ for such acceleration (Fig.~\ref{all}).

\begin{figure}[!b]
\centering\includegraphics[scale=0.9]{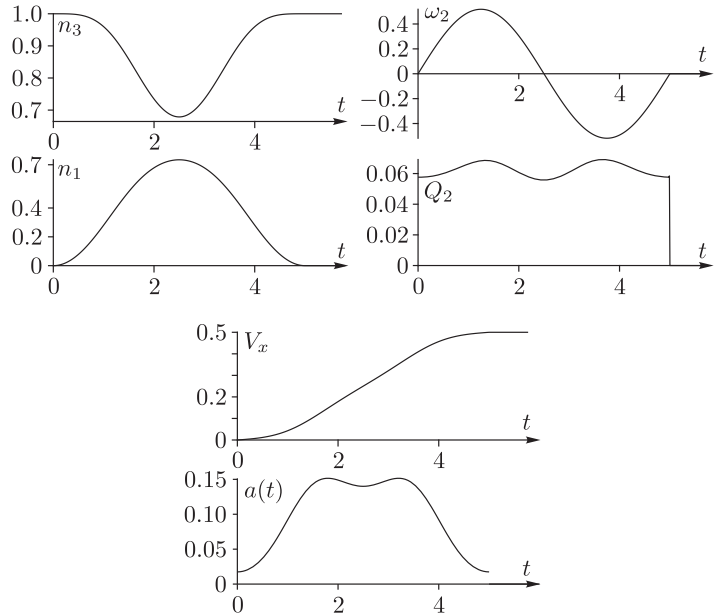}
\caption{\label{all}Time dependence of the vectors $\bn$, $\bs\omega$, $\bs Q$,
as well as the velocity and acceleration of the ball by acceleration for
one oscillation of the pendulum.}
\end{figure}

Figure~\ref{all} also shows the time-dependence of non-zero components of
the vector~$\bn$, velocity~$\bs V$ and acceleration $a(t)$ of the ball and
angular velocity of the pendulum~$\bom$. As seen in the figure, the ball
has gained the speed $\Delta V=0.5$ for the time interval $T=5$  moving
further with  constant velocity. The pendulum has executed one full
oscillation and returned to the initial condition.

To stop the ball, it is necessary to cause the pendulum to perform an
oscillation in the opposite direction with the same amplitude and for the
same time interval (i.e.~$\alpha=-0.83$ and~$T=5$). At the end of the
maneuver the pendulum is in the lower vertical position as during
acceleration.
\end{exl}

\subsection{Motion with change of the direction}

\begin{wrapfigure}[13]{o}{58mm}
\vspace{-3mm}
\centering\includegraphics{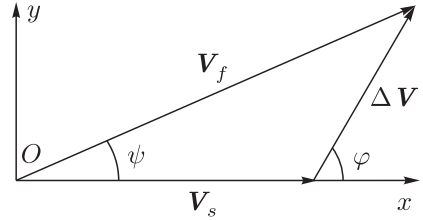}
\parbox{50mm}{\caption{\label{angle}$\bV_s$ is the initial velocity of  motion,
$\Delta \bV$ is the change of the velocity, $\bV_f$ is the velocity of the
ball after termination of a~maneuver.}}
\end{wrapfigure}
For simplicity we assume that at the initial instant of time the ball
moves along the axis $Ox$ with some constant velocity~$V_s$. We consider
a~motion where the pendulum executes one oscillation in an arbitrarily
predetermined direction (at an angle $\varphi$ to the direction of
motion, see Fig.~\ref{fig_theta}). The ball deviates from the initial
trajectory by some angle $\psi$ (see Fig.~\ref{angle}). The vector $\bn$
directed along the symmetry axis of the pendulum can be written as\vspace{0.5mm}
\[
\bn=\left(
\sin\theta\cos\varphi,
\sin\theta\sin\varphi,
\cos\theta
\right)^{\rm T},\vspace{0.5mm}
\]
where $\theta$ is the angle of deviation of the pendulum from the vertical, $\varphi$ is the angle between the initial direction of motion of the ball and the direction of  oscillation of the pendulum.

The equations of motion~\eqref{eqmQ} in projections onto the axes of the fixed coordinate system can be written as
\[
\begin{aligned}
&\dot{\Omega}_1I_0-R_oR_t(\dot{\omega}_1\cos\theta-\omega_1\sin\theta\dot{\theta})+Q_1=0,\\[1.5mm]
&\dot{\Omega}_2I_0-R_oR_t(\dot{\omega}_2\cos\theta-\omega_2\sin\theta\dot{\theta})+Q_2=0,\\[1.5mm]
&i_0\dot{\omega}_1-R_oR_t\dot{\Omega}_1\cos\theta-R_t\sin\theta\sin\varphi-Q_1=0,\\[1.5mm]
&i_0\dot{\omega}_2-R_oR_t\dot{\Omega}_2\cos\theta+R_t\sin\theta\cos\varphi-Q_2=0,\\[1.5mm]
&\dot{\theta}=\omega_2\cos\varphi-\omega_1\sin\varphi.
\end{aligned}
\]

The derivatives of the angular velocities of the ball  have the form
\beq
\dot{\Omega}_1=\dfrac{\dot{V}_2}{R_o}=\dfrac{a_2(t)}{R_o},\q
\dot{\Omega}_2=-\dfrac{\dot{V}_1}{R_o}=-\dfrac{a_1(t)}{R_o},\q
\dot{\Omega}_3=0.
\label{dOm12}
\eeq

We can represent the accelerations  $a_1(t)$ and $a_2(t)$ as
\beq
a_1(t)
=a(t,\alpha,T)\cos\varphi,\quad
a_2(t)=a(t,\alpha,T)\sin\varphi ,
\label{a12}
\eeq
where $a(t,\alpha,T)$ is defined by the expression~\eqref{acc} and $\theta$ is the function of  time and parameters $\alpha$ and $T$ and is defined by the expression~\eqref{theta} as before.

Differentiating Eq.~\eqref{rad-vector} with respect to time using~\eqref{dOm12} and~\eqref{a12} yields
\[
{\ddot{\bs r}}=R_o\bk\times{\dot{\bs\Omega}}=
R_o\left(-\dot{\Omega}_2,\,\dot{\Omega}_1,\,0\right)^{\rm T}=
a(t)\bs s,
\label{ddr_free}\]
where $\bs s=(\cos\varphi,\,\sin\varphi,\,0)^{\rm T}=\frac{\Delta\bV}{\Delta V}$ is the  unit vector (constant for an one maneuver) along which the velocity changes (see Fig.~\ref{angle}). Since $\varphi$ is the angle between the initial direction of motion and the direction of  oscillation of the pendulum, the following proposition holds.
\begin{pro}The velocity of the ball changes in the direction of  oscillation of the pendulum.\end{pro}

\begin{wrapfigure}[11]{o}{42mm}
\vspace{-4mm}
\centering\includegraphics[scale=0.9]{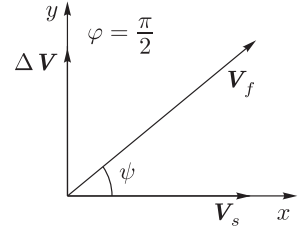}
\parbox{35mm}{\caption{\label{povorot_dv}Change of the velocity of the ball.}}
\end{wrapfigure}
If $\varphi=0$, we obtain acceleration of the ball in a~straight line considered in Section~\ref{razgon_po_pryamoi}.

Now we consider another case --- oscillation at the angle $\varphi=\frac{\pi}{2}$ to the initial direction of motion. Then the vector directed  along the symmetry axis of the pendulum can be represented as
\[
\bn=\left(
0,\,
\sin\theta,\,
\cos\theta
\right)^{\rm T}.
\]

The equations of motion~\eqref{eqmQ} can be rewritten as
\beq
\begin{aligned}
&I_0\dot{\Omega}_1-R_oR_t(\dot{\omega}_1\cos\theta+\omega_1^2\sin\theta)+Q_1=0,\\[1mm]
&i_0\dot{\omega}_1-R_oR_t\dot{\Omega}_1\cos\theta-R_t\sin\theta-Q_1=0,\\[1mm]
&\dot{\theta}=-\omega_1.
\end{aligned}
\label{eq_acc2}
\eeq

As in the previous case, we specify the angle of deviation of the pendulum from the
vertical in the form~\eqref{theta}. The derivatives of the angular velocities of the ball  have the
form\nopagebreak\vspace{1mm}
\beq
\dot{\Omega}_1=\dfrac{\dot{V}_2}{R_o}=\dfrac{a(t)}{R_o},\quad \dot{\Omega}_2=\dot{\Omega}_3=0.
\label{dOm1}\vspace{1mm}
\eeq
The expression for the acceleration of the ball $a(t,\alpha,T)$ has the same form as during acceleration in a~straight line~\eqref{acc} and the controls can be written as
\beq
\begin{gathered}
Q_1=-\dfrac{a(t,\alpha,T)}{R_o}\,I_0- R_o R_t\left(\cos\theta\ddot{\theta}-\sin\theta\dot{\theta}^2\right)\!, 
\\[1mm]
Q_2=Q_3=0.
\end{gathered}
\label{Q1}
\eeq
Thus, similarly to acceleration in a~straight line, for rotation through the specified angle it is necessary to generate a~control torque which is orthogonal to the plane of oscillation of the pendulum and, accordingly, to the direction of the vector of change of the velocity.

Differentiating Eq.~\eqref{rad-vector} with respect to time and using~\eqref{acc} and~\eqref{dOm1}, we obtain
\beq
{\ddot{\bs r}}=R_o\bk\times{\dot{\bs\Omega}}=R_o\left(-\dot{\Omega}_2,\,\dot{\Omega}_1,\,0\right)^{\rm T}=\left(0,\,a(t,\alpha,T),\,0\right)^{\rm T},
\label{ddr}\eeq
i.e., the velocity changes only in the direction of the axis $Oy$ --- in the direction of oscillation of the pendulum.

Our goal is to find a~value of the parameter $\alpha$ at which a~rotation through a~given angle $\psi$ is performed. To this end we integrate Eq.~\eqref{ddr} between $0$ and $T$ and compute the projections of the velocity $V_1$ and $V_2$ onto the axes $Ox$ and $Oy$, respectively, at the final instant of time, which are related with the angle of rotation in the absolute coordinate system by the relation (see Fig.~\ref{povorot_dv})
\[
\psi=\arctg\dfrac{V_{2}}{V_{1}}.
\]

\begin{figure}[!ht]
\begin{tabular}{cc}\small
\includegraphics{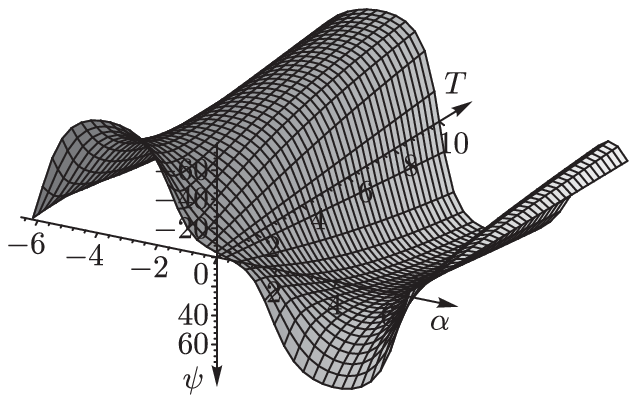}&
\includegraphics{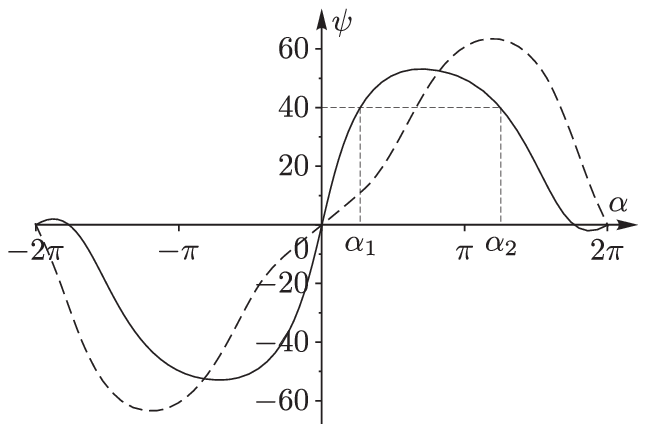}\\
({\it a})&({\it b})
\end{tabular}
\caption{\label{dphi}({\it a}) The surface $ \psi(\alpha,T)$. ({\it b}) Section formed
by the intersection of the surface $\psi(\alpha,T)$ with the planes $T=5$
(solid line) and $T=1$ (dashed line).}
\end{figure}

Changing the parameters $\alpha$ and $T$, we obtain the dependence $\psi(\alpha,T)$ shown in Fig.~\ref{dphi}({\it a}). The function $\psi(\alpha,T)$ is also antisymmetric relative to the plane $\alpha=0$, and $ \psi_{\max}\rightarrow\frac{\pi}{2}$ as $T\rightarrow0$.

Figure~\ref{dphi}({\it b}) shows sections of this surface for different values of $T$. As seen from the graphs,  at least two values of $\alpha$ correspond to the given $\psi$ and $T$ as in the previous case.

Substituting the computed value of the parameter $\alpha$, for example $\alpha_1$, into the equations of motion~\eqref{eq_acc2} and numerically integrating them using~\eqref{acc}, we obtain the controls~\eqref{Q1}.

\begin{figure}[!ht]
\centering\includegraphics[scale=0.9]{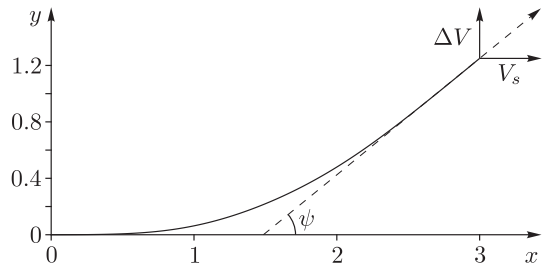}
\caption{\label{povorot}The trajectory of motion by rotation through the angle $\psi=40^{\circ}$, $V_s=0.6$.}
\end{figure}

\begin{exl}
We consider the rotation of the ball through the angle $\psi=40^{\circ}$ for the time $T=5$ with the initial conditions in the form
\[
\bn=\left(0,0,1\right)^{\rm T},\quad
\bs\omega=\left(0,0,0\right)^{\rm T},\quad
\bs r=\left(0,0,0\right)^{\rm T},\quad
\bs V=\left(V_s,0,0\right)^{\rm T},
\]
where $V_s$ is the initial velocity.

We choose the initial velocity such that the change of velocity $\Delta V$ and the parameter $\alpha$ are the same as in the previous example, i.e., $V_s =\Delta V\ctg\psi=0.6$.

The trajectory of such a~motion is shown in Fig.~\ref{povorot}.

Since all parameters are similar to those in the previous example, all
functions have the form shown in Fig.~\ref{all} up to the change of
variables $n_1\rightarrow n_2$, $\omega_2\rightarrow-\omega_1$,
$V_x\rightarrow V_y$ and $Q_2\rightarrow-Q_1$.
\end{exl}

\section{Motion with fixed parameters}\label{fix}

In this section we will consider a~control which prescribes the motion
with constant acceleration along a~given straight line (or
a~given curve in the more general case). Such a~control is used, for
example, in usual vehicles, i.e., in fact it is defined in a~body-fixed
reference frame. As will be shown below, in this case such an approach has
some shortcomings which stem from the conservatism of the resulting
system.


We demonstrate by a~specific example an algorithm for determining the
control torque. For this, we consider the uniformly accelerated motion of
the ball in a~straight line along the axis~$Ox$ under the law
\begin{equation*}
x(t)=\dfrac{a_0t^2}{2},
\end{equation*}
where $a_0=\const$ is the given acceleration of the ball.

Assume that the ball rolls without spinning, i.e., $\Omega_3=0$. Similarly to~\eqref{dOm2} using~\eqref{rad-vector} we obtain
\begin{equation*}
\dot{\Omega}_1(t)=0, \quad \dot{\Omega}_2(t)=-\dfrac{a_0}{R_o},
\end{equation*}
i.e., during acceleration along the axis $Ox$ the pendulum performs oscillations only in the plane~$Oxz$, therefore, the vector $\bn$ has the form
\[
\bn=
\left(
\sin\theta,\,
0,\,
\cos\theta
\right)^{\rm T},
\]
where $\theta$ is the angle of deviation of the pendulum from the
vertical, which is the unknown function of time. The angular velocity of
the pendulum has the form ${\omega_2=\dot{\theta}}$.

Substituting the resulting expressions into the equations of motion, we find the control torques which ensure uniformly accelerated motion in a~straight line along the axis $Ox$
\begin{gather*}
Q_1=Q_3=0,\\
Q_2=\dfrac{a_0}{R_o}I_0+R_oR_t\left(\cos\theta\ddot{\theta}-\sin\theta\dot{\theta}^2\right)
\end{gather*}
and the equation for determining the dependence of the angle $\theta$:
\beq
\ddot{\theta}=\dfrac{a_0I_0-a_0R_oR_t\cos\theta-R_oR_t\left(R_o\sin\theta\dot{\theta}^2+\sin\theta\right)}{R_o(i_0-R_oR_t\cos\theta)}.
\label{ddtheta}
\eeq

In addition to the geometrical integral $\bn^2=1$, this system admits another integral of motion quadratic in angular velocity
\[
\begin{aligned}
C&=-\frac{R_o}{2}\left(i_0-R_oR_t\cos\theta\right)^2\!\dot{\theta}^2+
R_oR_t\!\left(i_0-\frac{R_oR_t}{2}\cos\theta\right)\!
\left(\cos\theta-a_0\sin\theta\right)+{}\\
&{}+a_0\theta\left(I_0i_0+\frac{R_o^2R_t^2}{2}\right)-a_0R_oR_tI_0\sin\theta,
\end{aligned}
\]
which is an analog of an integral of generalized energy in a~uniformly accelerated reference frame.

\begin{figure}[!htb]
\centering\includegraphics[scale=0.9]{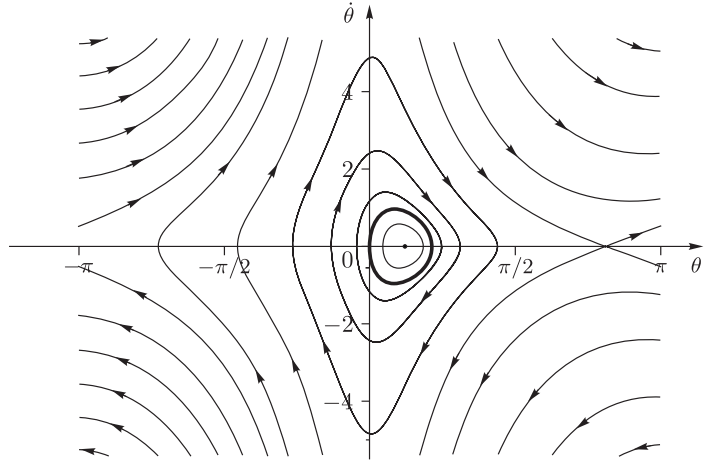}
\caption{\label{phase_razgon}Phase portrait of the system for uniformly accelerated motion of the ball in a~straight line.}
\end{figure}

Figure~\ref{phase_razgon} shows a~phase portrait of the
system~\eqref{ddtheta} on the plane $(\theta,\dot{\theta})$ for uniformly
accelerated motion of the ball in a~straight line with the acceleration
${a_0=0.1}$. For the given system parameters, there are two fixed points on
the phase plane one of which corresponds to a~stable equilibrium position
(at $\theta=0.38$, $\dot{\theta}=0$ of center type), while the other is an
unstable saddle point (at $\theta=2.56$ and $\dot{\theta}=0$).

\begin{figure}[!ht]
\centering\includegraphics[scale=0.9]{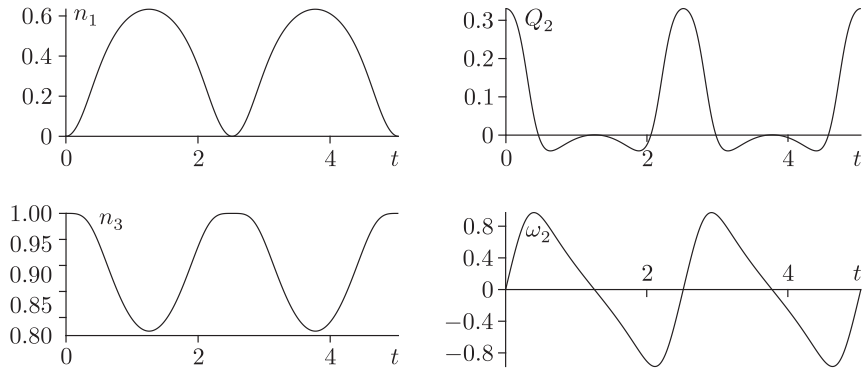}
\caption{\label{razgon}Time dependence of non-zero components of the
vectors $\bn,\bs\omega,\bs Q$ when the ball moves in a~straight line with
the constant acceleration $a_0=0.1$ and the initial conditions
$\theta=\dot{\theta}=0$.}
\end{figure}

As seen from Fig.~\ref{phase_razgon}, there exists a~single trajectory corresponding to a~state of rest at the initial instant of time (it passes through the point $(0,0)$ --- bold line). Periodic motion along this trajectory is ensured by the periodic control torque. Figure~\ref{razgon} shows the time dependence of the vectors $\bn,\bs\omega$ and $\bs Q$ for the initial conditions $\theta(0)=\dot{\theta}(0)=0$. From these graphs we can see that the pendulum performs oscillations in the plane $Oxz$, and the vector of the control torque changes periodically and is directed along the axis $Oy$.

If at the initial instant of time we define the angle of deviation  corresponding to the stable fixed point, then the ball will move with uniform acceleration. The deviation of the pendulum by a~constant angle will be maintained by a~constant value of the control torque. To any other (arbitrary) initial conditions there correspond closed periodic trajectories whose realization requires constantly applying the periodical control, which is inconvenient for the user.
That is, it is difficult to manually maintain the motion with  constant acceleration.

In addition, a~significant shortcoming of this control method is the  difficulty of switching to another motion mode at any instant of time (for example, if the velocity reaches the given value), because it may happen that the pendulum is not in the lower position.

The authors are grateful to A.\,V.\,Borisov, I.\,S.\,Mamaev and A.\,A.\,Kilin for fruitful discussions and useful remarks.

\end{document}